# Group Testing Problems
# in Experimental Molecular Biology

Preliminary Report

E. Knill
Los Alamos National Laboratory
knill@lanl.gov

S. Muthukrishnan
DIMACS
muthu@dimacs.rutgers.edu

March 1995

# 1  Introduction

In group testing, the task is to determine the *distinguished* members of a set of objects $\mathcal{L}$ by asking subset queries of the form "does the set $Q \subseteq \mathcal{L}$ contain a distinguished object?"

The primary biological application of group testing is for screening libraries of clones with hybridization probes [1, 2]. This is a crucial step in constructing physical maps and for finding genes. Group testing has also been considered for sequencing by hybridization [12]. Another important application includes screening libraries of reagents for useful chemically active zones.

In the literature on biological screening, the subset $Q$ which is used in a query is called a *pool*. The process of asking the query is referred to as *testing* or *screening* the pool with a probe. "Yes" answers and "No" answers are *positive* and *negative* results respectively. Accordingly, distinguished objects are called *positive* objects. The set of positive objects is denoted by $L_+$. Write $n = |\mathcal{L}|$.

The applications of group testing to biological screening have several features which determine the cost of finding the positive objects.

1. The same set of objects is screened with many different probes. Each probe is associated with its own set of positive objects.

2. It is expensive to prepare a pool for testing the first time, although once a pool is prepared, it can be screened many times with different probes.

3. Screening one pool at a time is expensive. Screening many pools in parallel with the same probe is much cheaper per pool.

4. It is common practice to individually test potential positive objects for confirmation. These confirmatory tests are relatively costly.

5. The screening results are not always reliable. Tests may be *false positive*, that is, they identify a positive in a pool when there does not exist any. Similarly, tests may be *false negative*, that is, they fail to identify a positive in a pool that contains positives. Therefore, errors must be tolerated.

6. In many cases it suffices to find just one or a few positive objects.

7. Prior information about the set of positives for a probe is best described in terms of a probability distribution on $2^{\mathcal{L}}$.

Features 1–3 strongly encourage nonadaptive approaches to screening. This implies that screening takes a few rounds where in each round as many pools as possible are screened in parallel. The cost of the strategy used is determined *both* by the number of pools that need to be constructed and by the number of tests that are performed (standard group testing problems, in contrast, consider only the latter cost.) Feature 4 arises when the chosen set of tests identify a superset of the positives (either due to errors in the tests or the particular choice of the set of tests); in this case, this superset is individually tested to isolate the set of positives. Feature 6 is not addressed in this survey, see [2]. Feature 7 entails average case analysis as opposed to the standard worst case analysis. In problems we state here, we will not necessarily add features 5 and 7 explicitly. Adding these features gives new problems that need to be investigated as well.

**Example of a screening effort.** One of the smaller screening efforts at the Human Genome Center at the Los Alamos National Laboratories involved finding positive clones in a library of 1298 clones ([3]). The expected number of positive clones for a probe was 2.6. 47 pools were prepared and they were screened with about 25 probes. It was discovered that there was substantial contamination compounding a false negative test rate estimated at 10% with a false positive test rate estimated at 7%. The effort of preparing the 47 pools required 10 hours of robot operation and several days for DNA extraction. About 10 screenings of the 47 pools can be accomplished in a day's work by an experienced technician. A similar number of individual confirmatory screenings of 5 to 10 clones per probe can be accomplished in a day.

In this particular screening effort, the error rates were higher than expected. Nevertheless, about 60% of the screenings yielded some positive clones with only a few individual confirmatory tests. The rest either had too many positive clones which required many confirmatory tests, or there was no detectable positive result. A Bayesean decoding technique was used, that is, the posterior probabilty of each clone's being positive was estimated by using suitable assumptions on the prior probability distribution of the positives. In this case it was sufficient to find one positive clone per probe, because different clones containing the same site were known to be identical for the purposes of this experiment with high probability. □

A number of novel group testing problems of relevance to molecular biology can be formulated using the features listed earlier. In what follows, we mention a few. We have chosen to present general versions of these problems; precise technical versions can be easily formulated. For further details, see [10].

## 2 General Problems: Non-adaptive Strategies

There are two main strategies for nonadaptive screening methods: the *strictly nonadaptive* [4] and the *trivial two-stage* strategy [8]. In the strictly nonadaptive strategy, a set of previously constructed pools is screened all at once. The positives must be determined from the results. The trivial two stage strategy is like the strictly nonadaptive strategy, except that candidate positive objects from screening the pools in the first round are individually confirmed in the second round; we call these the *confirmatory* tests. *Excess* confirmatory tests are those that involve negative objects. In practice, confirmatory tests are almost always used. The problems to be stated assume that the trivial two stage strategy is used.

In abstract group testing problems, it is assumed that the positive set $L_+$ satisfies some condition $A$. The best studied conditions are:

A1. The positive set is of cardinality exactly $p$.

A2. The positive set is of cardinality at most $p$.

If the group testing problem is probabilistic and the analysis of the screening method is average case, then uniform distributions over the set of allowed $L_+$ are assumed (unless stated otherwise).

**Problem 1** *What is the minimum number of pools required so that $L_+$ satisfying $A$ can always be determined with no more than $h$ excess confirmatory tests?*

**Problem 2** *Fix the number of pools to be $v$. Maximize $n$ such that for any $L_+$ satisfying $A$, $L_+$ can always be determined with no more than $h$ excess confirmatory tests.*

Both these are well-known problems (when $h = 0$) and background on these can be found in [4]. Good solutions to these problems have strong connections to finding optimal codes in coding theory. For example, the latter problem with assumption A2 is essentially that of constructing good superimposed codes or $p$-cover free families [7, 11]. There are instances of both these problems for which optimal bounds are not known.

The maximum $n$ in Problem 2 is denoted by $N(v, A, h)$ (no more than $h$ confirmatory tests in the worst case) and $\overline{N}(v, A, h)$ (no more than $h$ confirmatory tests on average). Bounds on these two quantities have been derived in [5, 6]. The most interesting asymptotics are those of

$$\overline{f}(p, h) = \limsup_{v \to \infty} \frac{\log_2 \overline{N}(v, A2), h)}{v} \qquad f(p, h) = \limsup_{v \to \infty} \frac{\log_2 N(v, A2, h)}{v}$$

**Problem 3** *Determine $f(p, h)$ and $\overline{f}(p, h)$.*

The current best bounds on $f(p, 0)$ obtained from [5, 6] are

$$\frac{\ln(2)}{p^2}(1 + o(1)) \leq f(p, 0) \leq 2\frac{\log_2(p)}{p^2}(1 + o(1)),$$

where the $o(1)$ is for $p \to \infty$. It has been shown that $\overline{f}(p, h)$ is independent of $h > 0$ [9] where the $o(1)$ is for $p \to \infty$. The best known bounds are

$$\frac{\ln(2)}{p}(1 + o(1)) \leq \overline{f}(p, 1) \leq \frac{1}{p}(1 + o(1))$$

The lower and upper bounds for $f$ (and $\overline{f}$) differ only by constants; however, note that this difference results in an exponential gap between the lower and upper bounds for $N$ (respectively, $\overline{N}$).

## 3 Pools vs Tests Tradeoff

Given the economics of screening clone libraries, there is a need to better understand the tradeoff between the number of pools that need to be constructed overall and the number of tests that need to be performed in the worst case. This tradeoff has not been investigated in the published group testing literature so far. Recently we have established some tradeoffs [10]. In what follows, we pose one technical problem in this context.

Consider *adaptive* group testing strategies, that is, those which perform a single test in each round; the outcome of a test in a round can determine the subsequent tests. Adaptive strategies can be naturally represented by a *decision tree*. The number of pools that need to be constructed to implement an adaptive strategy is the total number of distinct subsets that appear as queries in the decision tree; the worst case number of tests is the depth of the tree.

**Problem 4** *Consider a decision tree for determining $L_+$ satisfying an assumption $A$. Fix the number of pools to be $v$ and determine the minimum depth $D(v)$ of such a decision tree. Alternately, fix the depth $d$ and determine the minimum number of pools $V(d)$ in any such decision tree.*

# 4 Realistic Constraints

Considerations based on our biological screening setting sometimes helps us add constraints to the standard group testing problems resulting in novel and realistic variants. We state two such variants:

**The bounded bandwidth assumption.** There is a known ordering $x_1, \ldots, x_n$ of the objects such that for any $L_+ \in \mathcal{L}$, $L_+$ is included in an interval $\{x_i, x_{i+1}, \ldots x_{i+w-1}\}$ of length $w$. This assumption, denoted $B(w)$, arises in pooling libraries of clones which have already been mapped. The parameter $w$, called the *bandwidth*, is roughly related to the coverage of the clones, that is, the number of clones that overlap a particular site in the genome.

**The deBrujin-graph assumption.** This assumption is slightly technical to state. There is a known bijection of the objects to the words of length $l$ of a $\Sigma$-letter alphabet, such that $L_+$ corresponds to the set of all substrings of length $l$ of a length $w$ word. Here presence of defective objects is interdependent because of the overlap among substrings of length $l$. This is called the deBrujin graph assumption after the graphs that play a crucial role in modeling the defective objects. This problem arises in sequencing by hybridization [13].

Preliminary results based on assumption $B(w)$ for the general group testing problems in Section 2 and the tradeoffs in Section 3 can be found in [10]. More details on different types of problems based on the deBrujin-graph assumption are in [13].

# 5 Multiple probes

Interesting group testing problems arise when several inter-related probes for the given set of objects are considered at once. This scenario arises in physical mapping by unique markers where a library of clones is screened for a large number of markers to determine the relative ordering of the clones and the markers. For this purpose, it suffices to determine the incidence relationship between the probes and the clones. Here we abstract the group testing problem of finding these incidences.

In an idealized physical mapping setting, we are given a set of objects that are arbitrary unit length subintervals of an interval $I$ of length $L$. We are also given a set of probes which correspond to points on $I$. The task is to determine the probe-object incidences, that is, to determine the *incidence matrix* $M$ where $M(i,j)$ is 1 if the $i$th probe lies in the $j$th object and is 0 otherwise; the probes and objects are assumed to be ordered arbitrarily. Before we formulate problems for this task, we state a structural property of probe-object incidences. A matrix has the *monotone consecutive ones property* if the ones in each column are consecutive and the row indices $f(i)$ and $l(i)$ of the first and last one in column $i$, respectively, are non-decreasing functions of $i$. It is straightforward to see that for the incidence matrix obtained in the idealized physical mapping setting there exists permutations of the probes and the objects such that the permuted probe-object incidence matrix has the monotone consecutive ones property.

In general, in order to determine the probe-object matrix, one could construct pools of both probes and objects and screen them against each other. Here let us consider only the case of pooling objects and screening with single probes.

**Problem 5** *Given m probes and n objects whose incidence matrix is known to have the consecutive monotone ones property under some permutation, efficiently determine the incidences by pooling objects and screening the pools with probes.*

Group testing problems stated thus far are all pertinent to Problem 5: for instance, what is the minimum number of pools needed by a trivial two-stage strategy with at most $h$ confirmatory tests of individual probe-object incidences?

Another revalent scenario is when the objects and probes are assumed to be randomly drawn from the underlying interval $I$.

# 6 Discussion

Group testing abounds in experimental molecular biology. Indeed these applications have revitalized group testing research that was previously motivated primarily by blood screening, testing industrial products and multi-access communication problems. In this survey we have listed a small sample of the most tangible problems in the molecular biology setting. Among the other important issues omitted in this survey are average case analysis (that is, determining the expected cost under suitable assumption on the probabilistic distribution of the positives and Bayesian decoding of the pooling results) and tolerance to experimental noise.